\documentclass[final]{elsarticle}
\usepackage[letterpaper, margin=2cm]{geometry}
\usepackage[pdfstartview=FitH,bookmarks=false]{hyperref}
\usepackage{url}
\usepackage{breakurl}
\usepackage{graphicx,color}
\usepackage{paralist}
\usepackage{amsmath, amsthm, amssymb}
\usepackage{multirow}
\usepackage{fancyhdr}

\newcommand{\td}{\widetilde D}
\newcommand{\bx}{{\bf x}}
\newcommand{\by}{{\bf y}}
\newcommand{\rn}{\mathbb{R}^n}
\newcommand{\zn}{\mathbb{Z}^n}
\newcommand{\rr}{\mathbb{R}}
\newcommand{\sn}{S_{2}^{n-1}}

\begin{document}

\title{Comments on ``On Approximating Euclidean Metrics by Weighted $t$-Cost Distances in Arbitrary Dimension''}
\author{M.\ Emre Celebi}
\address{Department of Computer Science\\Louisiana State University, Shreveport, LA, USA\\
        \href{mailto:ecelebi@lsus.edu}{ecelebi@lsus.edu}}
\author{Hassan A.\ Kingravi}
\address{School of Electrical and Computer Engineering\\Georgia Institute of Technology, Atlanta, GA, USA\\
         \href{mailto:kingravi@gatech.edu}{kingravi@gatech.edu}}
         \author{Fatih Celiker}
\address{Department of Mathematics\\Wayne State University, Detroit, MI, USA\\
         \href{mailto:celiker@math.wayne.edu}{celiker@math.wayne.edu}}

\begin{abstract}

Mukherjee (Pattern Recognition Letters, vol. 32, pp. 824--831, 2011) recently introduced a class of distance functions called weighted $t$-cost distances that generalize $m$-neighbor, octagonal, and $t$-cost distances. He proved that weighted $t$-cost distances form a family of metrics and derived an approximation for the Euclidean norm in $\zn$. In this note we compare this approximation to two previously proposed Euclidean norm approximations and demonstrate that the empirical average errors given by Mukherjee are significantly optimistic in $\rn$. We also propose a simple normalization scheme that improves the accuracy of his approximation substantially with respect to both average and maximum relative errors.

\end{abstract}

\maketitle

\section{Introduction}
\label{sec_intro}

The Minkowski ($\mathcal{L}_p$) metric is inarguably one of the most commonly used quantitative distance (dissimilarity) measures in scientific and engineering applications. The Minkowski distance between two vectors $\bx = (x_1, x_2, \ldots, x_n)$ and {$\by = (y_1, y_2 , \ldots, y_n)$} in the $n$-dimensional Euclidean space, $\rn$, is given by

\begin{equation}
\label{equ_lp}
 \mathcal{L}_p (\bx,\by) =
 \left( {\sum\nolimits_{i = 1}^n {\left| {x_i  - y_i } \right| ^p } } \right)^{1/p}.
\end{equation}

Three special cases of the $\mathcal{L}_p$ metric are of particular interest, namely,
$\mathcal{L}_1$ (city-block metric), $\mathcal{L}_2$ (Euclidean metric),
and $\mathcal{L}_\infty$ (chessboard metric). Given the general form
(\ref{equ_lp}), $\mathcal{L}_1$ and $\mathcal{L}_2$ can be defined in a straightforward fashion, while $\mathcal{L}_\infty$ is defined as

\begin{equation*}
\mathcal{L}_{\infty}({\bf x},{\bf y}) = \max_{1\le i\le n}|x_i-y_i|.
\end{equation*}

The Minkowski metric enjoys the property of being translation invariant, i.e., $\mathcal{L}_p(\bx,\by) = \mathcal{L}_p(\bx + \bf{z},\by + \bf{z})$ for all $\bx, \by, \bf{z} \in \rn$. Since in many applications the data space is Euclidean, the most natural choice of metric is $\mathcal{L}_2$, which has the added advantage being isotropic (rotation invariant). For example, when the input vectors stem from an isotropic vector field, e.g., a velocity field, the most appropriate choice is to use the $\mathcal{L}_2$ metric so that all vectors are processed in the same way, regardless of their orientation \cite{Barni95}. However, $\mathcal{L}_2$ has the drawback of a high computational cost due to the multiplication and square root operations. As a result, $\mathcal{L}_1$ and $\mathcal{L}_\infty$ are often used as alternatives. Although these metrics are computationally more efficient, they deviate from $\mathcal{L}_2$ significantly.

Due to the translation invariance of $\mathcal{L}_p$, it suffices to consider $D_p(\bx) = \mathcal{L}_p (\bx,{\bf 0})$, i.e., the distance from the point $\bx$ to the origin. Therefore, in the rest of the paper, we will consider approximations to $D_p(\bx)$ rather than $\mathcal{L}_p (\bx,\by)$.

Let $\td$, defined on $\rn$, be an approximation to $D_2$ (Euclidean norm). We assume that $\td$ is a continuous and absolutely homogeneous function. Recall that $\td$ is called absolutely homogeneous (of degree one) if $\td(\lambda \bx) = |\lambda| \td(\bx) \quad \forall \lambda \in \rr, \; \forall \bx \in \rr^n.$

We note that all variants of $\td$ we consider in this paper satisfy these assumptions. As a measure of the quality of the approximation of $\td$ to $D_2$ we define the maximum relative error (MRE) as

\begin{equation}
\label{equ_max_err_temp1}
\varepsilon_{\text{max}}^{\td} = \sup \limits_{\bx \in \rn\setminus\{{\bf 0}\}}
{\frac{|{\td(\bx)-D_2(\bx)}|}{{D_2(\bx)}}}.
\end{equation}

Using the homogeneity of $D_2$ and $\td$, \eqref{equ_max_err_temp1} can be written as

\begin{equation}
\label{equ_max_err_temp2}
\varepsilon_{\text{max}}^{\td} = \sup_{\bx\in\sn} |\td(\bx)-1|,
\end{equation}

\noindent where $ \sn = \{ \bx \in \rn : D_2(\bx)=1 \} $ is the unit hypersphere of $\rn$ with respect to the Euclidean norm. Furthermore, by the continuity of $\td$, we can replace the supremum with maximum in \eqref{equ_max_err_temp2} and write

\begin{equation}
\label{equ_max_err}
\varepsilon_{\text{max}}^{\td} = \max_{\bx\in\sn} |\td(\bx)-1|.
\end{equation}

We will use \eqref{equ_max_err} as the definition of MRE throughout.

Mukherjee \cite{Mukherjee11} recently introduced a class of distance functions called weighted $t$-cost distances that generalize $m$-neighbor \cite{Das87a},  octagonal \cite{Das87b}, and $t$-cost \cite{Das92} distances. He proved that weighted $t$-cost distances form a family of metrics and derived an approximation for the Euclidean norm in $\zn$. Here we briefly review the $t$-cost norm.

The $t$-cost norm \cite{Das92} defines two points in the rectangular grid as neighbors when their respective hypercubes (or hypervoxels) share a hyperplane of any dimension. The cost associated with these points can be at most $t$, $1 \leq t \leq n$, such that if two consecutive points on a shortest path share a hyperplane of dimension $r$, the distance between them is taken as $\min(t, n - r)$. There are $n$ distinct $t$-cost norms defined by

\begin{equation*}
{D_t}(\bx) = \sum\limits_{i = 1}^t {{x_{(i)}}} ,\quad 1 \le t \le n
\end{equation*}

\noindent where $x_{(i)}$ is the $i$-th absolute largest component of $\bx$, i.e., $(x_{(1)}, x_{(2)}, \cdots, x_{(n)})$ is a permutation of $(|x_1|,|x_2|,\cdots,|x_n|)$ such that $x_{(1)} \geq x_{(2)} \geq \ldots \geq x_{(n)}$. The MRE of this norm is given by \cite{Das92}

\begin{equation*}
\varepsilon _{\max }^{{D_t}} = \max \left( {\sqrt {t} - 1 ,1 - \frac{t}{{\sqrt n }}} \right).
\end{equation*}

Mukherjee generalized the $t$-cost norm as follows \cite{Mukherjee11}:

\begin{equation*}
D_M(\bx) = \mathop {\max }\limits_{1 \le t \le n} \left\{ {{w_t}{D_t}(\bx)} \right\},
\end{equation*}

\noindent where $w_t$'s are non-negative real constants. Based on this weighted norm, he then derived an approximation for $D_2$ using the following weight assignment: ${w_t} = {1 \mathord{\left/ {\vphantom {1 {\sqrt t}}} \right. \kern-\nulldelimiterspace} {\sqrt t }}\;$ for $1 \leq t \leq n$. Note that $D_M$ consistently underestimates $D_2$ and the corresponding MRE is given by \cite{Mukherjee11}

\begin{equation}
\label{equ_mukherjee_mre}
\varepsilon _{\max }^{{D_{M}}} = 1 - \frac{1}{{\sqrt {\sum\nolimits_{i = 1}^n {{{\left( {\sqrt i  - \sqrt {i - 1} } \right)}^2}} } }}.
\end{equation}

In a recent study \cite{Celebi11}, we examined various Euclidean norm approximations in detail and compared their average and maximum errors using numerical simulations. Here we show that two of those approximations, namely Barni \emph{et al.}'s norm \cite{Barni95, Barni00} and Seol and Cheun's norm \cite{Seol08}, are viable alternatives to $D_M$.

Barni \emph{et al.} \cite{Barni95, Barni00} formulated a generic approximation for $D_2$ as

\begin{equation*}
D_B(\bx) = \delta \sum\limits_{i = 1}^n {\alpha _i x_{(i)}},
\end{equation*}

\noindent where $\boldsymbol{\alpha} = (\alpha_1, \alpha_2, \cdots, \alpha_n)$ and $\delta > 0$ are approximation parameters. Note that a non-increasing ordering and strict positivity of the component weights, i.e., $\alpha_1 \ge \alpha_2 \ge \cdots \ge \alpha_n > 0$ is a necessary and sufficient condition for $D_B$ to define a norm \cite{Barni00}.

Barni \emph{et al.} showed that the minimization of (\ref{equ_max_err}) is equivalent to determining the weight vector $\boldsymbol{\alpha}$ and the scale factor $\delta$ that solve the following minimax problem:

\begin{equation*}
\min\limits_{\boldsymbol{\alpha},\delta} \max\limits_{\bx \in V} \left| {D_B(\bx) - 1} \right|,
\end{equation*}

\noindent where $V = \{ \bx \in \rn \,:\, x_1 \ge x_2 \ge \cdots \ge x_n \ge 0,\; D_2(\bx)=1 \}$. The optimal solution and its MRE are given by

\begin{equation}
\label{equ_barni_mre}
\alpha_i^* = \sqrt{i}-\sqrt{i-1}, \qquad
\delta^* = \frac{2}{{1+\sqrt{\sum_{i=1}^n {{\alpha_i^*}^2}}}}, \qquad
\varepsilon_{\text{max}}^{D_B}= 1 - \delta^*.
\end{equation}

Note the striking similarity between \eqref{equ_mukherjee_mre} and \eqref{equ_barni_mre}. Interestingly, a similar but less rigorous approach had been published earlier by Ohashi \cite{Ohashi94}. It should also be noted that several authors approached the problem from a Euclidean distance transform perspective and derived similar approximations for the $2$- and $3$-dimensional cases, see for example \cite{Borgefors86} and \cite{Verwer91}. Furthermore, computation of weighted (Chamfer) distances in arbitrary dimensions on general point lattices is discussed in \cite{Fouard07}.

More recently, Seol and Cheun \cite{Seol08} proposed an approximation of the form

\begin{equation}
\label{equ_seol_cheun}
D_{a,b} (\bx) = a D_\infty(\bx) + b D_1(\bx),
\end{equation}

\noindent where $a$ and $b$ are strictly positive parameters to be determined by solving the following $2\times 2$ linear system

\begin{alignat*}{3}
 &a{\rm E}({D_\infty^2}  ) &&+ b{\rm E}({D_\infty D_1}) &&= {\rm E}({D_2 D_\infty}), \\
 &a{\rm E}({D_\infty D_1}) &&+ b{\rm E}({D_1^2}       ) &&= {\rm E}({D_2 D_1}),
\end{alignat*}

\noindent where $\rm E (\cdot)$ is the expectation operator.

Seol and Cheun estimated the optimal values of $a$ and $b$ using $100,000$ $n$-dimensional vectors whose components are independent and identically distributed, standard Gaussian random variables. In \cite{Celebi11}, we demonstrated that a fixed number of samples from the unit hypersphere gives biased estimates for the MRE. The basic reason behind this is the fact that a fixed number of samples fail to suffice as the dimension of the space increases.

It is easy to see that $D_B$ and $D_{a,b}$ fit into the general form

\begin{equation*}
\td(\bx) = \sum\limits_{i=1}^n {w_i x_{(i)}},
\end{equation*}

\noindent which is a weighted $D_1$ norm. For $D_B$ the weights are $w_1 = \delta^*$ and $w_{i \neq 1} = \delta^* \alpha_i^*$, whereas for $D_{a,b}$ they are
$w_1 = a + b$ and $w_{i \neq 1} = b$. Clearly, $D_B$ has a more elaborate design in which each component is assigned a weight proportional to its ranking (absolute value). However, this weighting scheme also presents a drawback in that a full ordering of the component absolute values is required.

$D_B$ and $D_{a,b}$ can also be written as linear combinations of the $D_1$ and $D_\infty$ norms, as in \eqref{equ_seol_cheun}. $D_1$ overestimates the $D_2$ norm, whereas $D_\infty$ underestimates it \cite{Chaudhuri92}. Therefore, it is natural to expect a suitable linear combination of $D_1$ and $D_\infty$ to give an approximation to $D_2$ better than either of them \cite{Rhodes95}. Note that Rosenfeld and Pfaltz \cite{Rosenfeld68} obtained a $2$-dimensional approximation by combining $D_1$ and $D_\infty$ nonlinearly as follows:
$\td(\bx) = \max \left( {\left\lfloor {{{2\left( {{D_1}(\bx) + 1} \right)} \mathord{\left/ {\vphantom {{2\left( {{D_1}(\bx) + 1} \right)} 3}} \right.
 \kern-\nulldelimiterspace} 3}} \right\rfloor ,{D_\infty(\bx)}} \right)$.

\section{Comparison of the Euclidean Norm Approximations}
\label{sec_comp}

\begin{figure}[!ht]
\centering
\includegraphics[width=0.6\columnwidth]{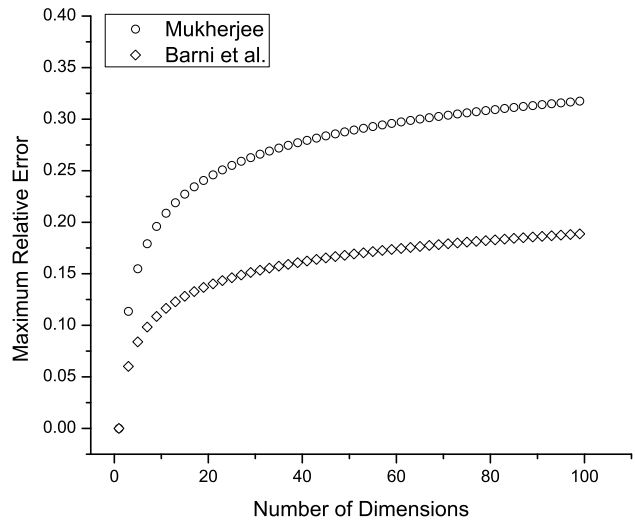}
\caption{\label{fig_max_err} Maximum relative errors for $D_M$ and $D_B$}
\end{figure}

Due to their formulations, the MREs for $D_M$ and $D_B$ can be calculated analytically using \eqref{equ_mukherjee_mre} and \eqref{equ_barni_mre}, respectively. In Figure \ref{fig_max_err} we plot the theoretical errors for these norms for $n \leq 100$. It can be seen that $D_B$ is not only more accurate than $D_M$, but also it scales significantly better.

The operation counts for each norm are given in Table \ref{tab_cost} (\textbf{ABS}: absolute value, \textbf{COMP}: comparison, \textbf{ADD}: addition, \textbf{MULT}: multiplication, \textbf{SQRT}: square root). The following conclusions can be drawn:

\begin{itemize}
\renewcommand{\labelitemi}{$\triangleright$}
    \item $D_B$ and $D_M$ have the highest computational cost due to the fact that they require sorting of the absolute values of the vector components.

    \item $D_{a,b}$ has the lowest computational cost among the approximate norms. A significant advantage of this norm is that it requires only two multiplications regardless of the value of $n$.

    \item $D_{a,b}$ can be used to approximate $D_2^2$ (squared Euclidean norm) using an extra multiplication. On the other hand, the computational cost of $D_B$ ($D_M$) is higher than that of $D_2^2$ due to the extra absolute value and sorting operations involved.
\end{itemize}

\begin{table}
\centering
\caption{ \label{tab_cost} Operation counts for the norms}
\vspace{.2cm}
\begin{tabular}{cccccc}
\hline
\textbf{Norm} & \textbf{ABS} & \textbf{COMP} & \textbf{ADD} & \textbf{MULT} & \textbf{SQRT}\\
\hline
\hline
$D_\infty$ & $n$ & $n - 1$ & 0 & 0 & 0\\
\hline
$D_1$ & $n$ & 0 & $n - 1$ & 0 & 0\\
\hline
$D_2$ & 0 & 0 & $n - 1$ & $n$ & 1\\
\hline
$D_B$ & $n$ & $\mathcal{O}(n \log{n})$ & $n - 1$ & $n$ & 0\\
\hline
$D_{a,b}$ & $n$ & $n - 1$ & $n$ & 2 & 0\\
\hline
$D_M$ & $n$ & $\mathcal{O}(n \log{n})$ & $n - 1$ & $n$ & 0\\
\hline
\end{tabular}
\end{table}

In Table \ref{tab_avg_max_err} we display the percentage average and maximum errors for $D_{a,b}$, $D_B$, and $D_M$ for $n \leq 8$. Average relative error (ARE) is defined as

\begin{equation*}
\varepsilon_{\text{avg}}^{\td} =
\frac{1}{|S|} \sum\limits_{\bx \in S}
{|\td(\bx) - 1|},
\end{equation*}

\noindent where $S$ is a finite subset of the unit hypersphere $S_2^{n-1}$, and $|S|$ denotes the number of elements in $S$. An efficient way to pick a random point on $S_2^{n-1}$ is to generate $n$ independent Gaussian random variables $x_1, x_2, \ldots, x_n$ with zero mean and unit variance. The distribution of the unit vectors

\begin{equation*}
\Big\{\by=(y_1,y_2,\ldots,y_n)\; : \; y_i={{x_i}/{\Big(\sum\nolimits_{j=1}^{n}x_j^2\Big)^{1/2},\quad i=1,2,\ldots,n}}\Big\}
\end{equation*}

\noindent will then be uniform over the surface of the hypersphere \cite{Muller59}. For each approximate norm, the ARE and MRE values were calculated over an increasing number of points, $ 2^{20}, 2^{21}, \ldots $ (that are uniformly distributed on the hypersphere) until the error values converge, i.e., the error values do not differ by more than $\epsilon = 10^{-5}$ in two consecutive iterations.

In Table \ref{tab_avg_max_err}, the error values under the column ``$D_{M}$ {\small{$(\rn)$}}" were obtained using the aforementioned iterative scheme, whereas those under the column ``$D_M$ {\small{$(\zn)$}}" are taken from \cite{Mukherjee11}. Motivated by the fact that $D_M$ consistently underestimates $D_2$, we also experimented with a normalized form of this approximate norm given by $D_{\widehat{M}}(\bx)=D_M(\bx) \mathord{\left/ {\vphantom {1 {{\delta^*}}}} \right. \kern-\nulldelimiterspace} {{\delta^*}}$. Note that $\delta^* < 1$ for $n \geq 2$ \eqref{equ_barni_mre}.

Note that for $D_M$ and $D_B$, two types of maximum error were considered: empirical maximum error ($\mbox{MRE}_e$), which is calculated numerically over $S$ and the theoretical maximum error ($\mbox{MRE}_t$), which is calculated analytically using \eqref{equ_mukherjee_mre} and \eqref{equ_barni_mre}, respectively.

\begin{table}
\centering
\small
{
\caption{ \label{tab_avg_max_err} Percentage average and maximum errors for the approximate Euclidean norms }
\begin{tabular}{|c|rr|rrr|rr|rrrrr|}
\hline
& \multicolumn{2}{|c|}{\multirow{2}{*}{$D_{a,b}$}} & \multicolumn{3}{|c|}{\multirow{2}{*}{$D_B$}} & \multicolumn{2}{|c|}{\multirow{2}{*}{$D_{\widehat{M}}$}} & \multicolumn{5}{|c|}{$D_{M}$}\\
\cline{9-13}
 & \multicolumn{2}{|c|}{} & \multicolumn{3}{|c|}{} & \multicolumn{2}{|c|}{} & \multicolumn{2}{|c}{{\small{$\rn$}}} & \multicolumn{2}{|c|}{{\small{$\zn$}}} &\\
\hline
\hline
 $n$ & ARE & $\mbox{MRE}_e$ & ARE & $\mbox{MRE}_e$ & $\mbox{MRE}_t$ & ARE & $\mbox{MRE}_e$ & ARE & $\mbox{MRE}_e$ & ARE & $\mbox{MRE}_e$ & $\mbox{MRE}_t$\\
\hline
2 & 2.00 & 5.25 & 2.41 & 3.96 & 3.96 & 2.48 & 4.12 & 2.55 & 7.61 & 2.40 & 7.61 & 7.61\\
\hline
3 & 2.39 & 9.98 & 3.00 & 6.02 & 6.02 & 2.97 & 6.40 & 4.14 & 11.35 & 3.63 & 11.35 & 11.35\\
\hline
4 & 2.57 & 13.64 & 3.44 & 7.39 & 7.39 & 3.28 & 7.97 & 5.21 & 13.75 & 4.29 & 13.75 & 13.75\\
\hline
5 & 2.68 & 16.59 & 3.77 & 8.39 & 8.39 & 3.53 & 9.16 & 5.98 & 15.47 & 4.65 & 15.46 & 15.49\\
\hline
6 & 2.73 & 18.88 & 4.01 & 9.19 & 9.19 & 3.73 & 10.12 & 6.55 & 16.80 & 4.85 & 16.79 & 16.83\\
\hline
7 & 2.76 & 20.67 & 4.18 & 9.84 & 9.84 & 3.92 & 10.91 & 7.00 & 17.90 & 5.00 & 17.86 & 17.92\\
\hline
8 & 2.77 & 21.92 & 4.31 & 10.39 & 10.39 & 4.10 & 11.59 & 7.35 & 18.78 & 5.04 & 18.75 & 18.82\\
\hline
\end{tabular}
}
\end{table}

By examining Table \ref{tab_avg_max_err}, the following observations can be made regarding the maximum error:

\begin{itemize}
\renewcommand{\labelitemi}{$\triangleright$}
    \item The most accurate approximation is $D_B$. This is because this norm is designed to minimize the maximum error.

    \item The proposed normalization is quite effective since the resulting norm, $D_{\widehat{M}}$, is, on the average, only $8.6$\% less accurate than $D_B$, whereas both $D_M$ {\small $(\rn)$} and $D_M$ {\small $(\zn)$} are, on the average, about $85$\% less accurate than $D_B$.

    \item The least accurate approximations are $D_M$ and $D_{a,b}$ for $n \leq 4$ and $n > 4$, respectively.

    \item As $n$ is increased, the error increases in all approximations. However, as can also be seen in Fig.\ \ref{fig_max_err}, the error grows faster in some approximations than others.

    \item For $D_B$, the empirical and theoretical errors agree almost perfectly in all cases, which demonstrates the validity of the presented iterative error calculation scheme. As for $D_M$, the agreement in each case is close, but not as close as that observed in $D_B$. We have confirmed that using a smaller convergence threshold ($\epsilon$) alleviates this problem at the expense of increased computational cost.
\end{itemize}

On the other hand, with respect to average error we can see that:

\begin{itemize}
\renewcommand{\labelitemi}{$\triangleright$}
    \item $D_{a,b}$ is the most accurate approximation. This is because this norm is designed to minimize the average error.

    \item $D_M$ {\small $(\rn)$} and $D_M$ {\small $(\zn)$} are the least accurate approximations. Furthermore, the errors given by Mukherjee are lower than those that we obtained (over $\rn$), and the discrepancy between the outcomes of the two error calculation schemes increases as $n$ is increased. The optimistic average error values given by Mukherjee are due to the fact that his approximation was primarily intended for use in digital geometry and hence the calculations were performed in $\zn$ (rather than $\rn$) using a very small number of points ranging from $32$ to $512$ \cite{Mukherjee11}. In fact, Mukherjee used progressively fewer points with increasing $n$ to calculate the error values. In \cite{Celebi11}, we demonstrated that more points are required in higher dimensions to obtain unbiased error estimates.
\end{itemize}

In the calculation of $D_{\widehat{M}}$, we assumed that the optimal scaling factor for $D_M$ is the same as that of $D_B$, i.e., $\delta^*$. In order to check this assumption, we performed a one-dimensional grid search over $\left[ \delta^*, 1\right]$ for each $n$ value. The results are shown in Table \ref{tab_avg_max_err_grid}. It can be seen that:

\begin{itemize}
\renewcommand{\labelitemi}{$\triangleright$}
    \item $D_{\widehat{M}}^{\widehat{\delta}}$ is significantly more accurate than $D_{\widehat{M}}^{\delta^*}$ with respect to both ARE and MRE.

    \item $D_{\widehat{M}}^{\widehat{\delta}}$ and $D_B$ have almost identical MREs. Since $D_B$ is analytically optimized for the maximum error it can be concluded that $D_{\widehat{M}}^{\widehat{\delta}}$ can reach the same optimality by means of a suitable scaling factor.

    \item Interestingly, $D_{\widehat{M}}^{\widehat{\delta}}$ is more accurate than $D_B$ with respect to ARE. This could be due to the fact that the two approximations take different paths towards minimizing the MRE.
\end{itemize}

\begin{table}
\centering
\small
{
\caption{ \label{tab_avg_max_err_grid} Percentage average and maximum errors for $D_{\widehat{M}}$}
\begin{tabular}{|c|rrr|rrr|}
\hline
& \multicolumn{3}{|c|}{\multirow{2}{*}{$D_{\widehat{M}}^{\delta^*}$}} & \multicolumn{3}{|c|}{\multirow{2}{*}{$D_{\widehat{M}}^{\widehat{\delta}}$}}\\
 & \multicolumn{3}{|c|}{} & \multicolumn{3}{|c|}{}\\
\hline
\hline
$n$ & ARE & $\mbox{MRE}_e$ & $\delta^*$ & ARE & $\mbox{MRE}_e$ & $\widehat{\delta}$\\
\hline
2 & 2.48 & 4.12 & 0.960434 & 2.41 & 3.96 & 0.961971\\
\hline
3 & 2.97 & 6.40 & 0.939809 & 2.79 & 6.02 & 0.943192\\
\hline
4 & 3.28 & 7.97 & 0.926150 & 2.99 & 7.39 & 0.931336\\
\hline
5 & 3.53 & 9.16 & 0.916059 & 3.13 & 8.40 & 0.922654\\
\hline
6 & 3.73 & 10.12 & 0.908117 & 3.23 & 9.18 & 0.915927\\
\hline
7 & 3.92 & 10.91 & 0.901603 & 3.31 & 9.84 & 0.910619\\
\hline
8 & 4.10 & 11.59 & 0.896101 & 3.40 & 10.39 & 0.905850\\
\hline
\end{tabular}
}
\end{table}

\section{Conclusions}
\label{sec_conc}

In this paper, we examined the weighted $t$-cost norm recently proposed by Mukherjee \cite{Mukherjee11} with respect to its ability to approximate the Euclidean norm in $\rn$. We evaluated the average and maximum errors of this norm using numerical simulations and compared the results to those of two other well-known Euclidean norm approximations. The results demonstrated that, because it was designed for digital geometry applications in $\zn$, the original weighted $t$-cost norm is not particularly suited to approximate the Euclidean norm in $\rn$. It is also shown, however, that when normalized with an appropriate scaling factor, Mukherjee's norm becomes competitive with an analytically optimized approximation with respect to both average and maximum relative errors.

\section{Acknowledgments}
\label{sec_ack}
This work was supported by grants from the Louisiana Board of Regents (LEQSF2008-11-RD-A-12) and US National Science Foundation (0959583, 1117457). The authors are grateful to the anonymous reviewers for their insightful suggestions and constructive comments that improved the quality and presentation of this paper.

\bibliographystyle{elsarticle-num}

\end{document}